\documentclass{amsart}[10pt]
\usepackage{amssymb}
\usepackage{amsfonts}
\usepackage{latexsym}
\usepackage{pb-diagram}

\newtheorem{theorem}{Theorem}[section]
\newtheorem{lemma}[theorem]{Lemma}
\newtheorem{proposition}[theorem]{Proposition}
\newtheorem{corollary}[theorem]{Corollary}
\theoremstyle{definition}
\newtheorem{definition}[theorem]{Definition}

\newtheorem{remark}[theorem]{Remark}

\newcommand{\cS}{{\mathcal S}}

\newcommand{\rank}{\text{rank\,}}

\newcommand{\Ad}{\operatorname{Ad}}

\newcommand{\Alt}{\operatorname{Alt}}

\newcommand{\nc}{\newcommand}
\nc{\Symm}{{\on{Sym}}}

\newcommand{\on}{\operatorname}

 \nc{\cE}{{\cal E}}

\renewcommand{\O}{{\mathcal{O}}}

\renewcommand{\a}{{\mathfrak a}}

\renewcommand{\c}{{\mathfrak c}}

\nc{\SL}{{\mathfrak sl}}
\nc{\HH}{{\mathfrak h}}
\newcommand{\g}{{\mathfrak{g}}}

\newcommand{\m}{{\mathfrak{m}}}

\nc{\wh}{\widehat}\nc{\wt}{\widetilde}

\newcommand{\ben}{\begin{enumerate}}
\newcommand{\een}{\end{enumerate}}

\newcommand{\cO}{{\mathcal O}}

\newcommand{\kk}{{\bf k}}

\newcommand{\AM}{{\mathbb{A}}}

\newcommand{\KM}{{\mathbb{K}}}

\renewcommand{\Ad}{{\operatorname{Ad}}}

\begin{document}

\title[Quantization of Poisson-Hopf stacks associated with group Lie bialgebras]
{Quantization of Poisson-Hopf stacks associated with group Lie
bialgebras}

\begin{abstract}
Let $G$ be a Poisson Lie group and $\g$ its Lie bialgebra. Suppose
that $\g$ is a group Lie bialgebra. This means that there is an
action of a discrete group $\Gamma$ on $G$ deforming the Poisson
structure into coboundary equivalent ones. Starting from this we
construct a non-trivial stack of Hopf-Poisson algebras and prove
the existence of associated deformation quantizations. This
non-trivial stack is a stack of functions on the formal Poisson
group, dual of the starting $\Gamma$ Poisson-Lie group. To
quantize this non-trivial stack we use quantization of a $\Gamma$
Lie bialgebra which is the infinitesimal of a $\Gamma$ Poisson-Lie
group (cf \cite{MS} for simple Lie groups and $\Gamma$ a covering
of the Weyl group and \cite{EH} for quantization in the general
case).
\end{abstract}

\author{Gilles Halbout}\thanks{\\
\noindent{MSC2000} Subject Classification Number: 17B37(Primary), 58H05(Secondary);
}
\address{Institut de Math\'ematiques et de Mod\'elisation de Montpellier,
Universit\'e de Montpellier 2, CC5149, Place Eug\`ene Bataillon, F-34095 Montpellier CEDEX 5, France}
\email{halbout@@math.univ-montp2.fr}

\author{Xiang Tang}
\address{Department of Mathematics, Washington University, St. Louis, Missouri, USA, 63130}
\email{xtang@@math.wustl.edu} \maketitle

\centerline{0. {\textsc{Introduction}}}

\bigskip

In this paper, we study examples of Poisson Hopf stacks and its
quantization. In \cite{EH}, the first author and his author
considered quantization of a $\Gamma$ Lie bialgebra ({\bf LBA}).
As an outcome of this quantization, they constructed a functor
from the category of $\Gamma$ Lie bialgebra to the category of
$\Gamma$ quantized universal enveloping algebras ({\bf QUE}). In
this paper, we first study the dual of a $\Gamma$ universal
enveloping algebra. Similar to the duality between Lie bialgebras
and Poisson-Lie groups, we discover a stack of Poisson formal
series Hopf algebras ({\bf PFSHA}), dual to a $\Gamma$ Lie
bialgebra. Then we study deformation quantization of this stack of
Poisson formal series Hopf algebras. We construct the deformation
quantization by applying the Drinfeld functor to a $\Gamma$
quantized universal enveloping algebra, and obtain a stack of
quantized formal series Hopf algebras ({\bf QFSHA}). We summarize
our results into the following commutative diagram.

\[
\begin{diagram}
\node{\Gamma\text{-\bf LBA}}\arrow{e,t}{\text{EH}}
\node{\Gamma\text{-\bf
QUE}}\arrow{s,r,<>}{\text{Dr}}\\
\node{\Gamma\text{-\bf
PFSHA}}\arrow{n,l,<>}{\approx}\arrow{e,t}{\text{Quant}}\node{\Gamma\text{-\bf
QFSHA}}
\end{diagram}
\]

\medskip

Let $\Gamma$ be a discrete group, $G$ a Lie group and $\g$ its Lie
algebra. Suppose that $\g$ is a $\Gamma$ Lie bialgebra (or
equivalently that $G$ is a $\Gamma$ Poisson group), i.e. a Lie
algebra $(\g,\mu_\g)$ together with a Lie cobracket $\delta_e$, an
action of $\Gamma$, $\theta : \Gamma \to \on{Aut}(\g,\mu_\g)$ and
$f : \Gamma \to \wedge^2(\g)$ a map satisfying compatibility rules
 such that $\Gamma$ acts on the double. Precise definitions and equivalent categories
corresponding to these objects will be recalled in Section 1.
Examples of $\Gamma$ Lie bialgebras arise from the following
situation: $G$ is a Poisson-Lie group with Lie bialgebra
$(\g,\mu_\g,\delta_\g)$, and $\Gamma \subset G$ is a discrete
subgroup. Another example is when $\g$ is a Kac-Moody Lie algebra
$\g$, and $\Gamma$ is a covering of the Weyl group of $\g$. In the
latter case, a quantization was given (\cite{MS}). Quantization of
a general $\Gamma$ Lie bialgebra was done in \cite{EH}. We will
also recall this quantization result in Section 1.

\medskip

It is then a natural question to ask what structure one gets on
the corresponding dual groups. Considering the function algebra of
a formal group, we get a trivial stack of Poisson Hopf algebras.
In the Section 3, we prove that we get a non-trivial stack of
Poisson algebras of functions on the formal Poisson Lie group
$G^*$ dual to a $\Gamma$ Poisson Lie group $G$. To do so, we will
construct ``lifts'' of the elements $(f(\gamma))_{\gamma \in
\Gamma}$ in the function algebra on $G^*$. In Section 2, we recall
basic definitions of stacks and of their quantizations.

\medskip

In Section 4, we construct a quantization of these non-trivial
Poisson-Hopf stacks. To do so we use quantization (cf \cite{EH})
of a $\Gamma$ Lie bialgebra. To deduce from it a quantization of a
non-trivial Poisson-Hopf stack we use the Drinfeld functor and
prove that quantization of the elements $(f(\gamma))_{\gamma \in
\Gamma}$ can be made ``admissible'' that is to say they will give
quantizations of the corresponding ``lifts''. Definitions of the
Drinfeld functor and admissibility will be recalled.

\medskip

Finally, in section 5, we give an explicit example corresponding
to the case where $G$ is a simple Lie group and $\Gamma$ a covering of the
corresponding Weyl group. In this case, quantization of Majid and
Soilbelman \cite{MS} will lead to an explicit quantization of the
non-trivial Poisson-Hopf stack.

\bigskip

Our results in this paper fit very well in the
Bressler-Gorokhovsky-Nest-Tsygan's framework \cite{bgnt} of
deformation quantization of gerbes. On one hand, our results
provide interesting examples of quantization of stacks, on the
other hand, the problems we are dealing with in this paper are
more special and complicated because we need to treat the Hopf
algebra structure. In \cite{KR} and \cite{SO} quantum Weyl groups
are used to study R-matrices, and we hope that the results in this
paper will shed a light on the general $\Gamma$ R-matrices.

\subsection*{Acknowledgements}

We thank Benjamin Enriquez for very useful discussions. The second
author would like to thank Institut de Math\'ematiques et de
Mod\'elisation de Montpellier for the hospitality of his visit. The research of the second author is partially supported by NSF grant 0604552.
\bigskip

\section{$\Gamma$ Lie bialgebras and equivalent categories}

In this section, we recall some results of \cite{EH}

\subsection{$\Gamma$ Lie algebras}

Define a group Lie algebra as a triple $(\Gamma,\g,\theta_\g)$,
where $\Gamma$ is a group, $\g$ is a Lie algebra and $\theta_\g :
\Gamma \to \on{Aut}(\g)$ is a group morphism. It is the
infinitesimal version of a $\Gamma$ action on a group $G$. Group Lie
algebras form a category.

If $\Gamma$ is a discrete group, a $\Gamma$ Lie algebra is a pair
$(\g,\theta_\g)$, such that $(\Gamma,\g,\theta_\g)$ is a group Lie
algebra. $\Gamma$ Lie algebras form a subcategory of group Lie
algebras. Such a $\Gamma$ Lie algebra will be said to be the
infinitesimal of a $\Gamma$ group $G$.

\medskip

Define a group cocommutative bialgebra as a triple $(\Gamma,U,i)$, where
$\Gamma$ is a group, $U$ is a cocommutative bialgebra,
$U = \oplus_{\gamma\in\Gamma} U_\gamma$ is a decomposition of $U$, and
$i : \kk\Gamma \to U$ is a bialgebra morphism, such that
$U_\gamma U_{\gamma'}\subset U_{\gamma\gamma'}$, $\Delta_U(U_{\gamma})
\subset U_\gamma^{\otimes 2}$, and $i$ is compatible with the $\Gamma$ grading.

We then define a $\Gamma$ cocommutative bialgebra as a pair $(U,i)$, such that
$(\Gamma,U,i)$ is a group cocommutative bialgebra. $\Gamma$ cocommutative
bialgebras form a
category.

The category of group (resp., $\Gamma$) cocommutative bialgebras contains
as a full subcategory the category of group (resp., $\Gamma$)
universal enveloping algebras, where $(U,\Gamma,i)$ satisfies the
additional requirement that $U_e$ is a universal enveloping algebra.

\medskip

Define a group commutative bialgebra (in a symmetric monoidal category $\cS$)
as a triple $(\Gamma,\cO,j)$, where $\Gamma$ is a group, $\cO$ is a
commutative algebra (in $\cS$) with a decomposition
$\cO = \oplus_{\gamma\in\Gamma} \cO_\gamma$,
such that $\cO_{\gamma}\cO_{\gamma'} = 0$ for $\gamma\neq\gamma'$,
algebra morphisms $\Delta_{\gamma'\gamma''} : \cO_{\gamma'\gamma''} \to
\cO_{\gamma'} \otimes \cO_{\gamma''}$,
$\eta : \kk \to \cO_e$ and $\varepsilon :\cO_e \to \kk$, satisfying
axioms such that when $\Gamma$ is finite, these morphisms add up to a
bialgebra structure on $\cO$; and $j : \cO \to \kk^\Gamma$ is a morphism
of commutative algebras, compatible with the $\Gamma$ gradings and the
maps $\Delta_{\gamma'\gamma''}$ on both sides. We define
$\Gamma$ commutative bialgebras as above.

We define the category of group (resp., $\Gamma$) formal series
Hopf (FSH) algebras as a full subcategory of the category of group
(resp., $\Gamma$) commutative bialgebras in $\cS = \{$pro-vector
spaces$\}$ by the condition the $\cO_e$ (or equivalently, each
$\cO_\gamma$) is a formal series algebra. Such FSH would
correspond to functions on the formal dual group of a $\Gamma$
group $G$.

\begin{proposition}\cite{EH}
1) We have (anti)equivalences of categories $\{$group Lie algebras$\}
\leftrightarrow \{$group universal enveloping algebras$\}
\leftrightarrow\{$group FHS algebras$\}$ (the last map is an antiequivalence).

2) If $\Gamma$ is a group, these (anti)equivalences restrict to
$\{\Gamma$-Lie algebras$\}
\leftrightarrow \{\Gamma$-universal enveloping algebras$\}
\leftrightarrow\{\Gamma$-FHS algebras$\}$.
\end{proposition}

If we denote the $\Gamma$ universal enveloping algebra
corresponding to a $\Gamma$ Lie algebra $(\Gamma,\g,\theta_\g)$ as
$U(\g) \rtimes \Gamma$. It is isomorphic to $U(\g) \otimes \kk\Gamma$
as a vector space;
if we denote by $x\mapsto [x]$, $\gamma\mapsto [\gamma]$ the natural maps
$\g\to U(\g) \rtimes \Gamma$, $\Gamma \to U(\g) \rtimes \Gamma$, then
the bialgebra structure of $U(\g) \rtimes \Gamma$ is given by
$[\gamma][x][\gamma^{-1}] = [\theta_\gamma(x)]$, $[\gamma][\gamma'] =
[\gamma\gamma']$, $[e] = 1$, $[x][x'] - [x'][x] = [[x,x']]$,
$\Delta([x]) = [x]\otimes 1 + 1 \otimes [x]$,
$\Delta([\gamma]) = [\gamma] \otimes [\gamma]$.

When $\Gamma$ is finite, the corresponding $\Gamma$ FSH algebra is then
$(U(\g) \rtimes \kk\Gamma)^*$, and in general, this is
$\oplus_{\gamma\in\Gamma} (U(\g) \otimes \kk\gamma)^*$.

\subsection{$\Gamma$ Lie bialgebras}

\begin{definition}
\label{defgamma}
A group Lie bialgebra
is a 5-uple $(\Gamma,\g, \theta_\g, \delta_\g,f)$
where $(\Gamma,\g,\theta_\g)$ is a group Lie algebra,
$\delta_\g : \g\to\wedge^2(\g)$ is\footnote{We view
$\wedge^2(V)$ as a  subspace of $V^{\otimes 2}$.} such that $(\g,\delta_\g)$
is a Lie bialgebra, and $f : \Gamma \to \wedge^2(\g)$
is a map $\gamma\mapsto f_\gamma$, such that:

a) $\wedge^2(\theta_\gamma) \circ \delta \circ
\theta_\gamma^{-1}(x) = \delta(x) + [f_\gamma,x\otimes 1 + 1 \otimes x]$
for any $x\in \g$,

b) $f_{\gamma\gamma'} = f_\gamma +
\wedge^2(\theta_\gamma)(f_{\gamma'})$,

c) $(\delta\otimes
\on{id})(f_\gamma)
+ [f_\gamma^{1,3},f_\gamma^{2,3}]$ + cyclic permutations  $= 0$.
\end{definition}
Group Lie bialgebras form a category. When $\Gamma$ is fixed, one defines
the category of $\Gamma$ Lie bialgebras as above.

\medskip

A co-Poisson structure on a group cocommutative bialgebra
$(\Gamma,U,i)$ is a co-Poisson structure $\delta_U : A \to \wedge^2(U)$,
such that $\delta_U(U_\gamma) \subset \wedge^2(U_\gamma)$.
Co-Poisson group cocommutative bialgebras form a category

Co-Poisson group universal enveloping algebras form a full subcategory of
the latter category. One defines the full subcategories of co-Poisson
$\Gamma$ cocommutative bialgebras and co-Poisson $\Gamma$ enveloping
algebras as above.

\medskip

A Poisson structure on a group commutative bialgebra
$(\Gamma,\cO,j)$ is a Poisson bialgebra structure $\{-,-\}:
\wedge^2(\cO) \to \cO$,
such that $\{\cO_\gamma,\cO_\gamma\} \subset \cO_{\gamma}$ and
$\{\cO_\gamma,\cO_{\gamma'}\} = 0$ if $\gamma\neq \gamma'$.
Poisson group bialgebras form a category, and Poisson group FSH
algebras form a full subcategory when $\cS = \{$pro-vector spaces$\}$.
One defines the full subcategories of Poisson $\Gamma$ bialgebras and
Poisson $\Gamma$ FSH algebras as above.

\medskip

{\bf Example.} Let $G$ be a Poisson-Lie (e.g., algebraic) group, let
$\Gamma \subset G$ be a subgroup (which we view as an abstract group).
We define $\theta_\gamma := \on{Ad}(\gamma)$, where $\on{Ad} : G \to
\on{Aut}_{Lie}(\g)$ is the adjoint action. If $P : G \to \wedge^2(\g)$ is the
Poisson bivector, satisfying $P(gg') = P(g') + \wedge^2(\on{Ad}(g'))(P(g))$,
then we set $f_\gamma := -P(\gamma)$. Then $(\g,\Gamma,f)$ is a
$\Gamma$ Lie bialgebra.

{\bf Example.} Assume that $(\g,r_\g)$ is a quasitriangular Lie
bialgebra and $\theta : \Gamma \to \on{Aut}(\g,t_\g)$ is an action
of $\Gamma$ on $\g$ by Lie algebra automorphisms preserving $t_\g
:= r_\g + r_\g^{2,1}$. If we set $f_\gamma :=
\theta_\gamma^{\otimes 2}(r) - r$, then $(\g,\theta,f)$ is a
$\Gamma$ Lie bialgebra (we call this a quasitriangular $\Gamma$
Lie bialgebra). For example, $\g$ is a Kac-Moody Lie algebra, and
$\Gamma = \widetilde{W}$ is a covering of the Weyl group of $\g$
(cf \cite{MS}).

\medskip

\begin{proposition}\cite{EH}
1) We have category (anti)equivalences $\{$group bialgebras$\}
\leftrightarrow \{$co-Poisson group universal enveloping
algebras$\}\leftrightarrow \{$Poisson group FSH algebras$\}$.

2) These restrict to category (anti)equivalences
$\{\Gamma$-bialgebras$\}
\leftrightarrow \{$co-Poisson $\Gamma$ universal enveloping
algebras$\}\leftrightarrow \{$Poisson $\Gamma$ FSH algebras$\}$.
\end{proposition}

If $(\g,\theta_\g,\delta_\g)$ is a $\Gamma$ Lie bialgebra,
then the co-Poisson structure on $U:= U(\g) \rtimes \Gamma$ is given by
$\delta_U([x]) = [\delta_\g(x)]$, and $\delta_U([\gamma])
= -[f_\gamma]([\gamma]\otimes [\gamma])$. (Here we also denote by
$x\mapsto [x]$ the natural map $\wedge^2(\g) \to \wedge^2(U(\g) \rtimes
\Gamma$).)

\subsection{Quantization of $\Gamma$ Lie bialgebras}

Let a $\Gamma$ graded bialgebra (in a symmetric monoidal category $\cS$)
be a bialgebra $A$ (in $\cS$), equipped with a grading
$A = \oplus_{\gamma\in \Gamma} A_\gamma$, such that $A_\gamma A_{\gamma'}
\subset A_{\gamma\gamma'}$ and $\Delta_A(A_{\gamma}) \subset
A_\gamma^{\otimes 2}$.

Assume that $A$ is a $\Gamma$ graded bialgebra in the category of
topologically free $\kk[[\hbar]]$-modules, quasicocommutative (in the sense
that $A_0 := A/\hbar A$ is cocommutative).
Then we get a co-Poisson structure on $A_0$. It is $\Gamma$ graded, in
the sense that $\delta_{A_0}((A_0)_\gamma) \subset \wedge^2((A_0)_\gamma)$.
We therefore get a classical limit functor $\on{class} :
\{\Gamma$-graded quasicocommutative
bialgebras$\} \to \{\Gamma$-graded co-Poisson bialgebras$\}$.

\begin{definition} A quantization functor for $\Gamma$ Lie bialgebras
is a functor $\{$co-Poisson $\Gamma$ universal enveloping
algebras$\} \to \{\Gamma$-graded quasicocommutative bialgebras$\}$,
right inverse to $\on{class}$.
\end{definition}

Assume that $(\g,\theta,f)$ is a $\Gamma$ Lie bialgebra. Let
$(U_e,*,\Delta_e)$ be the (Etingof-Kazhdan) quantization of ($\g$,
$\delta$)(we will also denote the multiplication by $m_e$). We get
from \cite{EH}:
\begin{proposition}
\label{propEH} There exist collections
$(\on{F}_{\gamma,\gamma\gamma'})_{\gamma,\gamma' \in \Gamma}$ of
elements in  $U^{\otimes 2}$ (with $\on{F}_{\gamma,\gamma\gamma'}
=1 + \hbar \on{f}_1 +O(\hbar^2)$ with
$\Alt(\on{f}_1)=\wedge^2(\theta_\gamma)(f_{\gamma'})$) \,
$(\on{v}_{\gamma,\gamma\gamma',\gamma\gamma'\gamma''})_{\gamma,\gamma',\gamma''
\in \Gamma}$ of elements in  $1 + \hbar^2U$,
$(U_\gamma,m_\gamma,\Delta_\gamma)_{\gamma \in \Gamma}$ of
bialgebras and $(\on{i}_{\gamma,\gamma\gamma'})_{\gamma,\gamma'
\in \Gamma}$ of algebra morphisms: $(U_\gamma, m_\gamma) \to
(U_{\gamma\gamma'}, m_{\gamma\gamma'})$ such that
\begin{itemize}
\item $\Delta_\gamma=\on{i}_{e,\gamma}^{\otimes 2}\circ
\on{Ad}(\on{F}_{e,\gamma})\circ \Delta_e \circ \on{i}_{e,\gamma}^{-1}$,
\item $(\on{F}_{e,\gamma} \otimes 1) *
(\Delta_e\otimes \on{id})(\on{F}_{e,\gamma}) =
(1\otimes \on{F}_{e,\gamma}) * (\on{id} \otimes \Delta_e)(\on{F}_{e,\gamma})$,
\item $\on{F}_{e,\gamma\gamma'} = \on{v}_{e,\gamma,\gamma\gamma'}^{\otimes 2} *
(\on{i}_{e,\gamma}^{\otimes 2})^{-1}(\on{F}_{\gamma,\gamma\gamma'}) *
\on{F}_{e,\gamma} * \Delta_e(\on{v}_{e,\gamma,\gamma\gamma'})^{-1}$,
\item $\on{i}_{e,\gamma\gamma'}= \on{i}_{\gamma,\gamma\gamma'} \circ \on{i}_{e,\gamma}
\circ \on{Ad}(\on{v}_{e,\gamma,\gamma\gamma'}^{-1})$,
\item $\on{v}_{e,\gamma\gamma',\gamma\gamma'\gamma''} * \on{v}_{e,\gamma,\gamma\gamma'} =
\on{v}_{e,\gamma,\gamma\gamma'\gamma''} * \on{i}_{e,\gamma}^{-1}(\on{v}_{\gamma,\gamma\gamma',\gamma\gamma'\gamma''})$.
\end{itemize}
Here $e$ is the unit of the group to make the formulas shorter but
could be any other element of the group and one would multiply
$\gamma$, $\gamma \gamma'$ and $\gamma\gamma'\gamma''$ on the left
by this elements in the formulas.
\end{proposition}
A quantization of the $\Gamma$ Lie bialgebra is then obtained as follows:
Set $U =
S(\g) \otimes \kk\Gamma[[\hbar]]$ and
$[x |\gamma] := x\otimes \gamma$, $[x\otimes x'|\gamma,\gamma']
:= (x\otimes \gamma) \otimes (x'\otimes \gamma')\in U^{\otimes 2}$.

There are unique linear maps $m : U^{\otimes 2}
\to U$ and $\Delta : U\to U^{\otimes 2}$, such that
$$
m : [x|\gamma][x'|\gamma'] \mapsto
[x * \on{i}_{e, \gamma}^{-1}(\theta_\gamma(x')) *
\on{v}_{e, \gamma, \gamma\gamma'}^{-1}
|\gamma\gamma']
$$
$$
\Delta : [x|\gamma] \mapsto [\Delta_e(x) *
\on{F}_{e, \gamma}^{-1}| \gamma,\gamma].
$$
The unit for $U$ is $[1|e]$, and the counit is the map
$[x|\gamma] \mapsto \delta_{\gamma,e} \varepsilon(x)$.

\begin{proposition}\cite{EH}
This defines a bialgebra structure on $U$,
quantizing the co-Poisson bialgebra structure induced by
$(\g,\theta,f)$.
\end{proposition}

\bigskip

\section{stack}

Let $M$ be a smooth manifold.
\begin{definition}
\label{defstack}
A stack on $M$ is the following data:
\begin{itemize}
\item an open cover of $M =\cup U_i$,
\item a sheaf of rings $A_i$ on every $U_i$,
\item an isomorphism of sheaves of rings $G_{ij}$: $A_j\vert (U_i\cap U_j)\to
A_i\vert (U_i \cap U_j)$
for every $i,j$,
\item  an invertible element $c_{ijk} \in A_i\vert(U_i \cap U_j \cap U_k)$ for every $i, j, k$
satisfying
\begin{itemize}
\item $G_{ij}G_{jk} = \on{Ad}(c_{ijk})G_{ik}$
\item and for every $i, j, k, l$,
$c_{ijk}c_{ikl} = G_{ij}(c_{jkl})c_{ijl}$.
\end{itemize}
\end{itemize}
If two such data $(U_i' , A_i', G_{ij}' , c_{ijk}')$ and $(U_i'' , A_i'', G_{ij}'' , c_{ijk}'')$ are
given on $M$, an isomorphism between them is
\begin{itemize}
\item an open cover $M =\cup U_i$
refining both $\{U_i'\}$ and $\{U_i'' \}$
\item isomorphisms $H_i$: $A_i'\to A_i''$
on $U_i$
\item and invertible elements $b_{ij}$ of $A_i'\vert (U_i \cap U_j)$ such that
\begin{itemize}
\item
$G_{ij}'' = H_i\on{Ad}(b_{ij})G_{ij}'H_j^{-1}$
\item and
$H_i^{-1}(c_{ijk}'') = b_{ij}G_{ij}'(b_{jk})c_{ijk}b_{ik}^{-1}$
\end{itemize}
\end{itemize}
\end{definition}

\bigskip

In what follows, we will still call a stack a collection of rings
$A_i$, group elements $G_{ij}$ and elements $c_{ijk}$ satisfying
the conditions above that is to say we will work without
considering the manifold $M$. More precisely, we will prove the
existence of a stack of Poisson Hopf algebra corresponding to
functions on the formal dual group $G^*$.
\begin{theorem}
\label{theoo1}
There exists a stack of Poisson Hopf algebras on $G^*$, i.e.:
\begin{itemize}
\item a collection $(\cO_{G^*_\gamma})_{\gamma \in \Gamma}$ of
Poisson Hopf algebras $(\cO_{G^*},m_0,\Delta_\gamma,\{-,-\}_\gamma)_{\gamma \in \Gamma}$,
\item Poisson morphisms $j_{\gamma,\gamma\gamma'}$: $\cO_{G^*_\gamma} \to \cO_{G^*_{\gamma\gamma'}}$,
\item elements $u_{\gamma,\gamma\gamma',\gamma\gamma'\gamma''}$ of $\cO_{G^*_\gamma}$
satisfying relations
\begin{itemize}
\item $j_{\gamma,\gamma\gamma'\gamma''}= j_{\gamma\gamma',\gamma\gamma'\gamma''} \circ j_{\gamma,\gamma\gamma'}
\circ
\on{Ad_{\star_\gamma}}(u_{\gamma,\gamma\gamma',\gamma\gamma'\gamma''}^{-1})$,
\item $u_{\gamma,\gamma\gamma'\gamma'',\gamma\gamma'\gamma''\gamma'''}
\star_\gamma u_{\gamma,\gamma\gamma',\gamma\gamma'\gamma''} =
u_{\gamma,\gamma\gamma',\gamma\gamma'\gamma''\gamma'''}
\star_\gamma
j_{\gamma,\gamma\gamma'}^{-1}(u_{\gamma\gamma',\gamma\gamma'\gamma'',\gamma\gamma'\gamma''\gamma'''})$.
\end{itemize}
\end{itemize}
The definition of the Baker-Campbell-Hausdorff product
$\star_\gamma$ will be recalled in
the next section.
\end{theorem}
Note that in this theorem (and the next one), one has to take
inverses of maps $j_{\gamma,\gamma\gamma'}$ and of elements
$u_{\gamma,\gamma\gamma',\gamma\gamma'\gamma''}$ to get equations
compatible with the ones of Definition \ref{defstack}

\smallskip

We will then prove the existence of a stack of algebras quantizing this stack
of Poisson Hopf algebras:
\begin{theorem}
\label{theoo2}
There exists a stack of algebras:
\begin{itemize}
\item $(\AM_\gamma,*_\gamma)_{\gamma \in \Gamma}$
quantizations\footnote{By quantization, we mean deformation
quantization, such that $\AM_\gamma/\hbar\AM_\gamma=
\cO_{G^*_\gamma}$, and $\frac{1}{\hbar}[\ ,\ ]_{*_\gamma}=\{\ ,\
\}_\gamma+O(\hbar)$.} of the Poisson algebras
$(\cO_{G^*_\gamma},\{-,-\}_\gamma )_{\gamma \in \Gamma}$,
\item algebra morphisms $i_{\gamma,\gamma\gamma'}$: $\AM_\gamma \to \AM_{\gamma\gamma'}$,
\item elements $v_{\gamma,\gamma\gamma',\gamma\gamma'\gamma''}$ of $\AM_\gamma$ such that elements $ev_{\gamma,\gamma\gamma',\gamma\gamma'\gamma''}:=\exp\left(\frac{v_{\gamma,\gamma\gamma',\gamma\gamma'\gamma''}}{\hbar}\right)$
satisfy relations
\begin{itemize}
\item $i_{\gamma,\gamma\gamma'\gamma''}= i_{\gamma\gamma',\gamma\gamma'\gamma''} \circ i_{\gamma,\gamma\gamma'}
\circ \on{Ad}(ev_{\gamma,\gamma\gamma',\gamma\gamma'\gamma''}^{-1})$,
\item $ev_{\gamma,\gamma\gamma'\gamma'',\gamma\gamma'\gamma''\gamma'''} *_\gamma ev_{\gamma,\gamma\gamma',\gamma\gamma'\gamma''} =
ev_{\gamma,\gamma\gamma',\gamma\gamma'\gamma''\gamma'''} *_\gamma i_{\gamma,\gamma\gamma'}^{-1}(ev_{\gamma\gamma',\gamma\gamma'\gamma'',\gamma\gamma'\gamma''\gamma'''})$.
\end{itemize}
\end{itemize}
\end{theorem}

\bigskip

\section{A Stack of Poisson bialgebras of functions on the formal group $G^*$}

Let $(\g, \theta_\g, \delta_\g,f)$ be a $\Gamma$ Lie bialgebra. In
this section we will construct a stack of Poisson bialgebras of
functions on a formal Poisson group $G^*$.

\subsection{Notations}

Let $(\g,\delta)$ be a Lie bialgebra and $(U(\g),\Delta_0,\delta)$
its corresponding cocommutative coPoisson bialgebra. The latter
can be seen as the dual of the function algebra of the formal
Poisson Lie group $G$ corresponding to $(\g,\delta)$. In the same
way, we will define $\cO_{G^*}$ as the commutative Poisson Hopf
algebra of functions of the formal Poisson Lie group $G^*$
corresponding to the dual Lie bialgebra $\g^*$.
We define by $\m_{G^*} \subset
\cO_{G^*}$ the maximal ideal of this ring. If $k$ is an integer $\geq 1$,
we denote by $\cO_{(G^*)^k}$ the
ring of formal functions on $(G^*)^k$, by $\m_{(G^*)^k}$ its
maximal ideal and by $\m_{(G^*)^k}^i$ the $i$-th power of this
ideal.

If $f,g\in \m_{(G^*)^k}^2$, then the series
$f \star g = f + g + {1\over 2} \{f,g\} +
\cdots + B_n(f,g) + \cdots$ is convergent, where $\sum_{i\geq 1}
B_i(x,y)$ is the Baker-Campbell-Hausdorff (BCH) series specialized to the
Poisson bracket of $\m_{(G^*)^k}^2$. The product $\star$ defines a  group structure on $\m_{(G^*)^k}^2$.

Let us recall a useful technical lemma (see \cite{EGH}, p. 2477),
proven for $m_{\g^*}$ and still true for $m_{G^*}$:
\begin{lemma}
\label{lemmetech}
For any $k\geq 1$ and $n \geq 2$,
$f,h \in \m_{(G^*)^k}^2$ and
$g \in \m_{(G^*)^k}^n$, one has
$$
f \star (h + g) = f \star h + g, \quad
(f+g)\star h = f \star h +g
\hbox{\ modulo\ }\m_{(G^*)^k}^{n+1}.
$$
\end{lemma}

\medskip

If $f\in \cO_{G^*}^{\wh\otimes n}$
and $P_1,\dots, P_m$ are disjoint subsets of
$\{1,\dots,m\}$, one defines $f^{P_1,\dots,P_n}$
using the coproduct of
$\cO_{G^*}$:
\begin{definition}
For $I_1,\dots,I_m$ disjoint ordered subsets of $\{1,\dots,n\}$,
$(U,\Delta)$ a Hopf algebra and $a \in  U^{\otimes m}$,
we define
$$a^{I_1,\dots,I_n}= \sigma_{I_1,\ldots,I_m} \circ
(\Delta^{|I_1|}\otimes \cdots \otimes \Delta^{|I_n|})(a),
$$
with $\Delta^{(1)}=\on{id}$, $\Delta^{(2)}=\Delta$,
$\Delta^{(n+1)}=({\on{id}}^{\otimes n-1} \otimes \Delta)\circ \Delta^{(n)}$,
and $\sigma_{I_1,\ldots,I_m} : U^{\otimes \sum_i |I_i|} \to
U^{\otimes n}$ is the
morphism corresponding to the map $\{1,\ldots,\sum_i |I_i|\}
\to \{1,\ldots,n\}$ taking $(1,\ldots,|I_1|)$ to $I_1$,
$(|I_1| + 1,\ldots,|I_1| + |I_2|)$ to $I_2$, etc.
\end{definition}
When $U$ is cocommutative, this definition depends only on the
underlying sets $I_1,\ldots,I_m$.

\medskip

When $(\g, \theta_\g, \delta_\g,f)$ is a $\Gamma$ Lie bialgebra we
thus get a collection of Lie bialgebras and so a collection
$(\cO_{G^*},m_0,\Delta_\gamma,\{-,-\}_\gamma)_{\gamma \in \Gamma}$
of Poisson bialgebras. We will denote by $\star_\gamma$ the
corresponding BCH products.

\medskip

\subsection{``Lifts'' and functional equations}

We will now construct ``lifts'' $\tilde{f}_{\gamma,\gamma\gamma'}
\in \m_{G^*}^{\wh\otimes 2}$ of the elements
$\wedge^2(\theta_\gamma)(f_{\gamma'})$, $\gamma,\gamma' \in
\Gamma$ that will satisfy similar relation as
$\on{F}_{\gamma,\gamma\gamma'}$ in Proposition \ref{propEH}.
\begin{proposition}
\label{lift}
Let $\gamma,\gamma'$ be in $\Gamma$. Then
there exists $\tilde{f}_{\gamma,\gamma\gamma'}$ in  $m_{G^*}^{\wh\otimes 2}$ the image of which in
$\g^{\otimes 2}$ under the square of
the projection $\m_{G^*} \to \m_{G^*} / \m_{G^*}^2 = \g$
equals $\wedge^2(\theta_\gamma)(f_{\gamma'})$, and such that
\begin{equation}
\label{twistequation} (\tilde{f}_{\gamma,\gamma\gamma'} \otimes 1)
\star_\gamma (\Delta_\gamma\otimes
\on{id})(\tilde{f}_{\gamma,\gamma\gamma'}) = (1\otimes
\tilde{f}_{\gamma,\gamma\gamma'}) \star_\gamma  (\on{id} \otimes
\Delta_\gamma)(\tilde{f}_{\gamma,\gamma\gamma'}).
\end{equation}
Such a $\tilde{f}_{\gamma,\gamma\gamma'}$ is unique up to the action of
$\m^2_{G^*}$
by $\lambda \cdot \tilde{f} = \lambda^{1} \star_\gamma \lambda^{2} \star_\gamma \tilde{f}
\star_\gamma (-\lambda)^{12}$.
We will call such a $\tilde{f}$ a twist for $\Delta_\gamma$.
\end{proposition}

{\em Proof.} Let us construct $\tilde{f}_{\gamma,\gamma\gamma'}$
by induction: we will construct a convergent sequence $\tilde{f}_N
\in \m_{G^*}^{\wh\otimes 2}$ ($N\geq 2$) satisfying
(\ref{twistequation}) in $\m_{G^*}^{\wh\otimes 3} /
(\m_{G^*}^{\wh\otimes 3} \cap \m_{(G^*)^3}^N)$, where $\m_{(G^*)^3}^N$ is the $N$-th power of $\m_{(G^*)^3}$. When $N = 3$, we take for $\tilde{f}_2$ any
lift of $\wedge^2(\theta_\gamma)(f_{\gamma'})$ to
$\m_{G^*}^{\wh\otimes 2}$; then equation (\ref{twistequation}) is
automatically satisfied.

To shorten the notation, we will write $\tilde{f}_{1,2}$ for
$\tilde{f}_{\gamma,\gamma\gamma'}$, $\tilde{f}_{2,3}$ for
$\tilde{f}_{\gamma\gamma',\gamma\gamma'\gamma''}$ and so on and
the same thing for $\alpha_{-,-,-}$

Let $N$ be an integer $\geq 3$; assume that we have constructed
$\tilde{f}_N$ in $\m_{G^*}^{\wh\otimes 2}$ satisfying equation
(\ref{twistequation}) in $\m_{G^*}^{\wh\otimes 3}/
(\m_{G^*}^{\wh\otimes 3}\cap\m_{(G^*)^3} ^N)$. Set $\alpha^N
_{1,2,3}:= \tilde{f}^N_{1,2}\star_\gamma \tilde{f}^N_{12,3}
-\tilde{f}^N_{2,3} \star_\gamma \tilde{f}^N_{1,23}$. Then
$\alpha^N_{1,2,3}$ belongs to $\m_{G^*}^{\wh\otimes
3}\cap\m_{(G^*)^3}^N$, and the following equalities hold in
$\m_{G^*}^{\wh\otimes 4}/(\m_{G^*}^{\wh\otimes 4} \cap
\m_{(G^*)^4}^{N+1})$:
\[
\begin{array}{rl}
\alpha^N_{12,3,4} = &~ \tilde{f}^N_{1,2} \star_\gamma
\alpha^N_{12,3,4}= \tilde{f}^N_{1,2} \star_\gamma
\tilde{f}^N_{12,3}\star_\gamma \tilde{f}^N_{123,4} -
\tilde{f}^N_{1,2} \star_\gamma \tilde{f}^N_{3,4} \star_\gamma
\tilde{f}^N_{12,34}\\
\\
=&~\alpha^N_{1,2,3}+ \tilde{f}^N_{2,3} \star_\gamma
\tilde{f}^N_{1,23} \star_\gamma
\tilde{f}^N_{123,4}-\tilde{f}^N_{3,4}\star_\gamma
\tilde{f}^N_{1,2}\star_\gamma \tilde{f}^N_{12,34}\cr
\\
& ~(\hbox{using Lemma \ref{lemmetech}})\\
\\
=&~\alpha^N_{1,2,3} + \tilde{f}^N_{2,3} \star_\gamma
\tilde{f}^N_{1,23} \star \tilde{f}^N_{123,4}
-\tilde{f}^N_{3,4}\star (\tilde{f}^N_{2,34}\star_\gamma
\tilde{f}^N_{1,234} +\alpha^N_{1,2,34}) \cr
\\
&~(\hbox{using Lemma \ref{lemmetech}}\hbox{ and the definition of
}\alpha^N_{1,2,34})\cr
\\
=&~\alpha^N_{1,2,3} + \tilde{f}^N_{2,3} \star_\gamma
(\alpha^N_{1,23,4}+ \tilde{f}^N_{23,4} \star_\gamma
\tilde{f}^N_{1,234})\cr
\\
&- \alpha^N_{1,2,34}-\tilde{f}^N_{3,4}\star_\gamma
\tilde{f}^N_{2,34}\star_\gamma \tilde{f}^N_{1,234} \cr
\\
&~(\hbox{using the definition of }\alpha^N_{1,23,4} \hbox{ and
Lemma \ref{lemmetech}})\cr
\\
=&~\alpha^N_{1,2,3} + \alpha^N_{1,23,4}+ (\tilde{f}^N_{3,4}
\star_\gamma \tilde{f}^N_{2,34}+ \alpha^N_{2,3,4} )\star_\gamma
\tilde{f}^N_{1,234} \cr
\\
&- \alpha^N_{1,2,34}-\tilde{f}^N_{3,4}\star_\gamma
\tilde{f}^N_{2,34}\star_\gamma \tilde{f}^N_{1,234} \cr
\\
&~
(\hbox{using the definition of }\alpha^N_{2,3,4} \hbox{ and Lemma
\ref{lemmetech}})\cr
\\
=&~  \alpha^N_{1,2,3} + \alpha^N_{1,23,4} - \alpha^N_{1,2,34} +
\alpha^N_{2,3,4}\cr
\\
&~ (\hbox{using Lemma \ref{lemmetech}}).
\end{array}
\]

Let us denote by $\overline\alpha^N$ the image of $\alpha^N$ in
$(\m_{\g^*}^{\wh\otimes 3}\cap\m_{(\g^*)^3}^N) /
(\m_{\g^*}^{\wh\otimes 3}\cap\m_{(\g^*)^3}^{N+1}) =
(S^{>0}(\g)^{\otimes 3})_N$, then we get
$$
\overline\alpha^N_{12,3,4} + \overline\alpha^N_{1,2,34}
=\overline\alpha^N_{1,2,3} + \overline\alpha^N_{1,23,4}
+\overline\alpha^N_{2,3,4}.
$$
This means that $\bar{\alpha}$ is a cocycle for the subcomplex
$(S^{>0}(\g)^{\otimes\cdot},d)$ of the co-Hochschild complex.
Using \cite{Dr:QH}, Proposition 3.11, one proves that the $k$-th
cohomology group of this subcomplex is $\wedge^k(\g)$, and that
the antisymmetrization map coincides with the canonical projection from
the space of cocycles to the cohomology group. For $N=3$, the
equations of Definition \ref{defgamma}
 implies
$\on{Alt}(\overline\alpha^3)=0$, and hence $\overline\alpha^3$ is
the coboundary of an element $\overline\beta_3\in
(S^{>0}(\g)^{\otimes 2})^3$. For $N>3$, $\overline\alpha^N$ is the
coboundary of an element $\overline\beta^N \in
(S^{>0}(\g)^{\otimes 2})^N$, since the degree $N$ part of the
cohomology vanishes. We then set $\tilde{f}^{N+1} := \tilde{f}^N +
\beta^N$, where $\beta^N\in\m_{G^*}^{\wh\otimes 2} \cap
\m_{(G^*)^2}^N$ is a representative of $\overline\beta^N$. Then
$\tilde{f}^{N+1}$ satisfies (\ref{twistequation}) in
$\m_{G*}^{\wh\otimes 3} /(\m_{G^*}^{\wh\otimes 3} \cap
\m_{(G^*)^3}^{N+1})$.

The sequence $(\tilde{f}_N)^{N\geq 2}$ has a limit $\tilde{f}$,
which then satisfies (\ref{twistequation}).

\medskip

The second part of the theorem can be proved in the same way or by analyzing the
choices for $\overline\beta_N$ in the above proof.
\hfill \qed \medskip

\subsection{Isomorphism of formal Poisson manifolds $G_\gamma^* \simeq G_{\gamma\gamma'}^*$}

\begin{proposition}
\label{iso} Let $\gamma,\gamma' \in \Gamma$ and let $G_\gamma^*$
and  $G_{\gamma\gamma'}^*$ be the formal Poisson-Lie groups associated to the
corresponding Lie cobrackets. There exists an isomorphism of
Poisson algebras $j_{\gamma,\gamma\gamma'}$: $\cO_{G_\gamma^*}
\simeq \cO_{G_{\gamma\gamma'}^*}$.
\end{proposition}

{\em Proof.} Let $P : \wedge^2(\cO_{G_\gamma^*}) \to \cO_{G_\gamma^*}$ be the Poisson
bracket on $\cO_{G_\gamma^*}$ corresponding to the
Lie-Poisson Poisson
structure on $G_\gamma^*$.
Then $(\O_{G_\gamma^*},m_0,P,\Delta_\gamma)$ is a Poisson formal series Hopf (PFSH)
algebra; it corresponds to the formal Poisson-Lie group $G_\gamma^*$
equipped with its Lie-Poisson structure.

Set ${}^{\tilde{f}_{\gamma,\gamma\gamma'}}\Delta_\gamma(a) =
\tilde{f}_{\gamma,\gamma\gamma'}\star_\gamma\Delta_\gamma(a)\star_\gamma
(-\tilde{f}_{\gamma,\gamma\gamma'})$ for any
$a\in\cO_{G_\gamma^*}$. It follows from the fact that
$\tilde{f}_{\gamma,\gamma\gamma'}$ satisfies the equation
(\ref{twistequation}) that
$(\cO_{G_\gamma^*},m_0,P,{}^{\tilde{f}_{\gamma,\gamma\gamma'}}\Delta_\gamma)$
is a PFSH algebra.

Let us denote by ${\bf PFSHA}$ and ${\bf LBA}$ the categories of PSFH algebras
and Lie bialgebras. We have a category equivalence
$c : {\bf PFSHA} \to {\bf LBA}$, taking $(\cO,m,P,\Delta)$ to the Lie bialgebra
$(\c,\mu,\delta)$, where $\c := \m/\m^2$ ($\m\subset\cO$ is the maximal ideal),
the Lie cobracket of $\c$ is induced by $\Delta - \Delta^{2,1} : \m\to
\wedge^2(\m)$, and the Lie bracket of $\c$ is induced by the Poisson
bracket $P : \wedge^2(\m) \to \m$. The inverse of the functor $c$
takes $(\c,\mu,\delta)$ to $\cO = \wh S(\c)$ equipped with its usual product;
$\Delta$ depends only on $\delta$ and $P$ depends on $(\mu,\delta)$.

Then $c$ restricts to a category equivalence $c_{\on{fd}} : {\bf
PFSHA}_{\on{fd}} \to {\bf LBA}_{\on{fd}}$ of subcategories of
finite-dimensional objects (in the case of ${\bf PFSH}$, we say
that $\cO$ is finite-dimensional if and only if $\m/\m^2$ is).

Let $\on{dual} : {\bf LBA}_{\on{fd}} \to {\bf LBA}_{\on{fd}}$ be
the duality functor. It is a category antiequivalence; we have
$\on{dual}(\g,\mu,\delta) = (\g^*,\delta^t,\mu^t)$. Then
$\on{dual} \circ c_{\on{fd}} : {\bf PFSHA}_{\on{fd}} \to {\bf
LBA}_{\on{fd}}$ is a category antiequivalence. Its inverse is the
usual functor $\g\mapsto U(\g)^*$. If $G$ is the formal
Poisson-Lie group with Lie bialgebra $\g$, one sets $\cO_G =
U(\g)^*$.

Let us apply the functor $c$ to
$(\cO_{G_\gamma^*},m_0,P,{}^{\tilde{f}_{\gamma,\gamma\gamma'}}\Delta_\gamma)$.
We obtain $\c = \m/\m^2 = \g$; the Lie bracket is unchanged with
respect to the case $\tilde{f}_{\gamma,\gamma\gamma'}=0$, so it is
the Lie bracket of $\g$; the Lie cobracket is given by
$\delta_{\gamma\gamma'}(x) =
\delta_{\gamma}+[\wedge^2(\theta_\gamma)(f_{\gamma'}),x\otimes 1 +
1\otimes x]$ since the reduction of
$\tilde{f}_{\gamma,\gamma\gamma'}$ modulo
$(\m_{G_\gamma^*})^2\wh\otimes \m_{G_\gamma^*} + \m_{G_\gamma^*}
\wh\otimes (\m_{G_\gamma^*})^2$ is equal to
$\wedge^2(\theta_\gamma)(f_{\gamma'})$.

Then applying $\on{dual} \circ c_{\on{fd}}$ to
$(\cO_{G_\gamma^*},m_0,P,{}^{\tilde{f}_{\gamma,\gamma\gamma'}}\Delta_\gamma)$,
we obtain the Lie bialgebra $(\g^*,\delta_{\gamma\gamma'})$. So
this PFSH algebra is isomorphic to the PFSH algebra of the formal
Poisson-Lie group $G^*_{\gamma\gamma'}$. Let us call this PFSH
algebra morphism $j_{\gamma,\gamma\gamma'}$.

In particular, the Poisson algebras
$\cO_{G_\gamma^*}$ and $\cO_{G_{\gamma\gamma'}^*}$ are isomorphic.
\hfill \qed \medskip

\begin{remark}
It is easy to check that the map $\g=\m_{G_\gamma^*}\slash \m^2_{G_\gamma^*}
\to \m_{G_{\gamma\gamma'}^*}\slash \m^2_{G_{\gamma\gamma'}^*}=\g$ induced
by the isomorphism $j_{\gamma,\gamma\gamma'}$ is the identity.
\end{remark}

\begin{remark}
We have proven a stronger result than the existence of a Poisson
algebra morphism $j_{\gamma,\gamma\gamma'}$: $\cO_{G_\gamma^*}
\simeq \cO_{G_{\gamma\gamma'}^*}$. This morphism intertwines the
coproducts in the following way:
$$\Delta_{\gamma\gamma'}=j_{\gamma,\gamma\gamma'}^{\otimes 2}
\circ{}^{\tilde{f}_{\gamma,\gamma\gamma'}}\Delta_\gamma\circ j_{\gamma,\gamma\gamma'}^{-1}.$$
\end{remark}

\subsection{Composition of equivalences}

Let us first prove the following lemma:
\begin{lemma}
For $\gamma,\gamma'$ in $\Gamma$, the element
$(j_{\gamma,\gamma\gamma'}^{\otimes
2})^{-1}(\tilde{f}_{\gamma\gamma',\gamma\gamma'\gamma''})\star_\gamma
\tilde{f}_{\gamma,\gamma\gamma'}$ is a solution of the equation
\begin{equation}
(\tilde{f} \otimes 1) \star_\gamma
(\Delta_\gamma\otimes \on{id})(\tilde{f}) =
(1\otimes \tilde{f}) \star_\gamma  (\on{id} \otimes \Delta_\gamma)(\tilde{f}).
\end{equation}
\end{lemma}

{\em Proof.} One can check that directly or notice that
$\tilde{f}_{\gamma\gamma',\gamma\gamma'\gamma''}$ is a twist for
$\Delta_{\gamma\gamma'}$. Therefore
$(j_{\gamma,\gamma\gamma'}^{\otimes
2})^{-1}(\tilde{f}_{\gamma\gamma',\gamma\gamma'\gamma''})$ is a
twist for $(j_{\gamma,\gamma\gamma'}^{\otimes 2})^{-1} \circ
\Delta_{\gamma\gamma'} \circ j_{\gamma,\gamma\gamma'}=
{}^{\tilde{f}_{\gamma,\gamma\gamma'}}\Delta_\gamma$. Accordingly
$(j_{\gamma,\gamma\gamma'}^{\otimes
2})^{-1}(\tilde{f}_{\gamma\gamma',\gamma\gamma'\gamma''})\star_\gamma
\tilde{f}_{\gamma,\gamma\gamma'}$ is a twist for $\Delta_\gamma$.
\hfill \qed \medskip

Let us then notice that the image of
$(j_{\gamma,\gamma\gamma'}^{\otimes
2})^{-1}(\tilde{f}_{\gamma\gamma',\gamma\gamma'\gamma''})\star_\gamma
\tilde{f}_{\gamma,\gamma\gamma'}$ under the square of the
projection $\m_{G^*} \to \m_{G^*} / \m_{G^*}^2 = \g$ equals
$\wedge^2(\theta_\gamma)(f_{\gamma'})+\wedge^2(\theta_{\gamma\gamma'})
(f_{\gamma''})=\wedge^2(\theta_\gamma)(f_{\gamma'}+\wedge^2(\theta_{\gamma'})(f_{\gamma''}))=
\wedge^2(\theta_\gamma)(f_{\gamma'\gamma''})$. Thanks to
Proposition \ref{lift}, there exists an element
$u_{\gamma,\gamma\gamma',\gamma\gamma'\gamma''}$ in $1+
\m_{G^*}^2$ such that
$$\tilde{f}_{\gamma,\gamma\gamma'\gamma''}=u_{\gamma,\gamma\gamma',\gamma\gamma'\gamma''}^{\otimes 2}
\star_\gamma(j_{\gamma,\gamma\gamma'}^{\otimes
2})^{-1}(\tilde{f}_{\gamma\gamma',\gamma\gamma'\gamma''})\star_\gamma
\tilde{f}_{\gamma,\gamma\gamma'}\star_\gamma
\Delta_\gamma(u_{\gamma,\gamma\gamma',\gamma\gamma'\gamma''})^{-1}.$$

Finally, from the previous section, we defined
$j_{\gamma,\gamma\gamma'} $,
$j_{\gamma\gamma',\gamma\gamma'\gamma''} $ and
$j_{\gamma,\gamma\gamma'\gamma''} $ such that
\begin{equation}
\begin{array}{rl}
\Delta_{\gamma\gamma'\gamma''}=&j_{\gamma,\gamma\gamma'\gamma''}^{\otimes
2}\circ{}^{\tilde{f}_{\gamma,\gamma\gamma'\gamma''}}\Delta_\gamma\circ
j_{\gamma,\gamma\gamma'\gamma''}^{-1}\cr
\\
=&j_{\gamma,\gamma\gamma'\gamma''}^{\otimes
2}\circ{}^{u_{\gamma,\gamma\gamma',\gamma\gamma'\gamma''}^{\otimes
2}\star_\gamma (j_{\gamma,\gamma\gamma'}^{\otimes
2})^{-1}(\tilde{f}_{\gamma\gamma',\gamma\gamma'\gamma''})\star_\gamma
\tilde{f}_{\gamma,\gamma\gamma'}\star_\gamma
\Delta_\gamma(u_{\gamma,\gamma\gamma',\gamma\gamma'\gamma''})^{-1}}\Delta_\gamma\circ
j_{\gamma,\gamma\gamma'\gamma''}^{-1}\cr
\\
=&(j_{\gamma,\gamma\gamma'\gamma''}\circ
\Ad_{\star_\gamma}(u_{\gamma,\gamma\gamma',\gamma\gamma'\gamma''}))^{\otimes
2}\circ{}^{ (j_{\gamma,\gamma\gamma'}^{\otimes
2})^{-1}(\tilde{f}_{\gamma\gamma',\gamma\gamma'\gamma''})\star_\gamma
\tilde{f}_{\gamma,\gamma\gamma'}}\Delta_\gamma\cr
\\
&\circ (j_{\gamma,\gamma\gamma'\gamma''}\circ
\Ad_{\star_\gamma}(u_{\gamma,\gamma\gamma',\gamma\gamma'\gamma''}))^{-1}\cr
\\
=&(j_{\gamma,\gamma\gamma'\gamma''}\circ
\Ad_{\star_\gamma}(u_{\gamma,\gamma\gamma',\gamma\gamma'\gamma''})\circ
j_{\gamma,\gamma\gamma'}^{-1}\circ
j_{\gamma\gamma',\gamma\gamma'\gamma''}^{-1})^{\otimes 2}\circ
\Delta_{\gamma\gamma'\gamma''}\cr
\\
& \circ (j_{\gamma,\gamma\gamma'\gamma''}\circ
\Ad_{\star_\gamma}(u_{\gamma,\gamma\gamma',\gamma\gamma'\gamma''})\circ
j_{\gamma,\gamma\gamma'}^{-1}\circ
j_{\gamma\gamma',\gamma\gamma'\gamma''}^{-1})^{-1}.
\end{array}
\end{equation}
By the equivalence $c_{\on{fd}}$ between the category {\bf PFSHA$_{\on{fd}}$} and {\bf LBA$_{\on{fd}}$}  we get
$$j_{\gamma,\gamma\gamma'\gamma''}=
j_{\gamma\gamma',\gamma\gamma'\gamma''}\circ
j_{\gamma,\gamma\gamma'}\circ
\Ad_{\star_\gamma}(u_{\gamma,\gamma\gamma',\gamma\gamma'\gamma''}^{-1}).$$

\subsection{Cocycle relation for the $u_{\gamma,\gamma\gamma',\gamma\gamma'\gamma''}$}

We will end this section by proving the following proposition that will prove
Theorem \ref{theoo1}:
\begin{proposition}
For any $\gamma,\gamma',\gamma'',\gamma'''$ in $\Gamma$, we have
$$u_{\gamma,\gamma\gamma'\gamma'',\gamma\gamma'\gamma''\gamma'''}
\star_\gamma u_{\gamma,\gamma\gamma',\gamma\gamma'\gamma''} =
u_{\gamma,\gamma\gamma',\gamma\gamma'\gamma''\gamma'''}
\star_\gamma
j_{\gamma,\gamma\gamma'}^{-1}(u_{\gamma\gamma',\gamma\gamma'\gamma'',\gamma\gamma'\gamma''\gamma'''}).$$
\end{proposition}
{\em Proof.} To shorten the notation, we will write
$\tilde{f}_{1,2}$ for $\tilde{f}_{\gamma,\gamma\gamma'}$,
$\tilde{f}_{2,3}$ for
$\tilde{f}_{\gamma\gamma',\gamma\gamma'\gamma''}$ and so on and
the same thing for the $j_{-,-}$ and the $u_{-,-,-}$. We will omit
the BCH product $\star_\gamma$ and write $\star$ for the product
$\star_{\gamma\gamma'}$, $\Delta_0$ for the coproduct
$\Delta_\gamma$ and $\Delta$ for the coproduct
$\Delta_{\gamma\gamma'}$. We will also write $j(-)$ instead of
$j^{\otimes 2}(-)$ when no confusion is possible.

\smallskip

We have by definition $\tilde{f}_{1,4}\Delta_0 u_{1,3,4}=u_{1,3,4}^{\otimes 2}
j_{1,3}^{-1}(\tilde{f}_{3,4})\tilde{f}_{1,3}$. Multiplying this equality
on the right by $\Delta_0 u_{1,2,3}$ and using the fact that
$\tilde{f}_{1,3}\Delta_0 u_{1,2,3}=u_{1,2,3}^{\otimes 2}
j_{1,2}^{-1}(\tilde{f}_{2,3})\tilde{f}_{1,2}$, we
 get
\[
\tilde{f}_{1,4}\Delta_0 u_{1,3,4}\Delta_0
u_{1,2,3}=u_{1,3,4}^{\otimes 2}
j_{1,3}^{-1}(\tilde{f}_{3,4})u_{1,2,3}^{\otimes 2}
j_{1,2}^{-1}(\tilde{f}_{2,3})\tilde{f}_{1,2}.
\]
Using now that $j_{1,3}^{-1}(-)u_{1,2,3}=u_{1,2,3}j_{1,2}^{-1}
\circ j_{2,3}^{-1}(-)$, we get
\begin{equation}
\label{equu1} \tilde{f}_{1,4}\Delta_0 u=u^{\otimes 2}j_{1,2}^{-1}
\circ j_{2,3}^{-1}(\tilde{f}_{3,4})
j_{1,2}^{-1}(\tilde{f}_{2,3})\tilde{f}_{1,2},
\end{equation}
where $u=u_{1,3,4}u_{1,2,3}$. On the other hand, we have
$\tilde{f}_{2,4}\star \Delta u_{2,3,4}=u_{2,3,4}^{\otimes 2}\star
j_{2,3}^{-1}(\tilde{f}_{3,4})\star \tilde{f}_{2,3}$. Using the
Poisson algebra morphism $j_{1,2}$ and the fact that
$j_{1,2}^{-1}\circ \Delta =\tilde{f}_{1,2}
\Delta_0(j_{1,2}^{-1}(-))\tilde{f}_{1,2}^{-1}$, we get
\begin{equation}
\label{equu2} j_{1,2}^{-1}(\tilde{f}_{2,4})
\tilde{f}_{1,2}\Delta_0(j_{1,2}^{-1}(u_{2,3,4}))\tilde{f}_{1,2}^{-1}
=j_{1,2}^{-1}(u_{2,3,4}^{\otimes 2}) j_{1,2}^{-1} \circ
j_{2,3}^{-1}(\tilde{f}_{3,4})j_{1,2}^{-1}( \tilde{f}_{2,3}).
\end{equation}
From $\tilde{f}_{1,4}\Delta_0u_{1,2,4}=u_{1,2,4}^{\otimes
2}j_{1,2}^{-1}(\tilde{f}_{2,4})\tilde{f}_{1,2}$, using Equation
(\ref{equu2}), we get
\begin{equation}
\label{equu3}
\tilde{f}_{1,4}\Delta_0(u')=
(u')^{\otimes 2}
j_{1,2}^{-1} \circ j_{2,3}^{-1}(\tilde{f}_{3,4})j_{1,2}^{-1}( \tilde{f}_{2,3})\tilde{f}_{1,2},
\end{equation}
where $u'=u_{1,2,4}j_{1,2}^{-1}(u_{2,3,4})$. Then Equations
(\ref{equu1}) and (\ref{equu3}) imply that if $w=u(u')^{-1}$ then
$\tilde{f}_{1,4}\Delta_0(w)=w\tilde{f}_{1,4}$, and so if
$w'=j_{1,4}(w)$ then $\Delta_0(w')=w'$. Recall that by similar
properties of $u_{i,j,k}$, $w' \in 1 +m_{G^*}^2$. Suppose that
$w'\not=1$ and set $i\geq 2$ the largest possible $i$ such that
$w' \in 1 +m_{G^*}^i$, but not in $1+m_{G^*}^{i+1}$. Let
$\bar{w}'$ be the projection of $w'$ in $m_{G^*}^i\slash
m_{G^*}^{i+1}$. Relation $\Delta_0(w')=w'$ implies that $\bar{w}'$
is in $\g$ and so in $m_{G^*}^1$ which is a contradiction. Thus we
have proved that $w=w'=1$ and so that $u=u'$.
\section{quantization}

\subsection{Duality of QUE and QFSH algebras} \label{sect:duality}

In this subsection, we recall some facts from \cite{Dr:QG} (proofs
can be found in \cite{Gav}). Let us denote by ${\bf QUE}$ the
category of quantized universal enveloping (QUE) algebras and by
${\bf QFSH}$ the  category of quantized formal series Hopf (QFSH)
algebras. We denote by ${\bf QUE}_{\on{fd}}$ and ${\bf
QFSH}_{\on{fd}}$ the subcategories corresponding to finite
dimensional Lie bialgebras.

We have contravariant functors ${\bf QUE}_{\on{fd}} \to {\bf QFSH}_{\on{fd}}$,
$U\mapsto U^*$ and ${\bf QFSH}_{\on{fd}} \to {\bf QUE}_{\on{fd}}$,
$\cO\mapsto \cO^\circ$. These functors are inverse to each other.
$U^*$ is the full topological dual of $U$, i.e., the space of all
continuous (for the $\hbar$-adic topology) $\KM[[\hbar]]$-linear maps
$U \to \KM[[\hbar]]$.
$\cO^\circ$ the space of continuous $\KM[[\hbar]]$-linear forms
$\cO\to \KM[[\hbar]]$,
where $\cO$ is equipped with the $\m$-adic topology (here $\m\subset \cO$
is the maximal ideal).

We also have covariant functors ${\bf QUE} \to {\bf QFSH}$, $U\mapsto U'$
and ${\bf QFSH} \to {\bf QUE}$, $\cO\mapsto \cO^\vee$. There functors are
also inverse to each other. $U'$ is a subalgebra of $U$, while $\cO^\vee$
is the $\hbar$-adic completion of $\sum_{k\geq 0} \hbar^{-k} \m^k \subset
\cO[1/\hbar]$.

We also have canonical isomorphisms $(U')^\circ \simeq (U^*)^\vee$
and $(\cO^\vee)^* \simeq (\cO^\circ)'$.

If $\a$ is a finite dimensional Lie bialgebra and $U = U_\hbar(\a)$
is a QUE algebra quantizing $\a$, then $U^* = \cO_{A,\hbar}$ is a
QFSH algebra quantizing the Poisson-Lie group $A$ (with Lie bialgebra $\a$),
and $U' = \cO_{A^*,\hbar}$ is a QFSH algebra quantizing the Poisson-Lie
group $A^*$ (with Lie bialgebra $\a^*$). If now $\cO = \cO_{A,\hbar}$
is a QFSH algebra quantizing $A$, then $\cO^\circ = U_\hbar(\a)$ is a
QUE algebra quantizing $\a$ and $\cO^\vee = U_\hbar(\a^*)$ is a QFSH algebra
quantizing $\a^*$.

We now compute these functors explicitly in the case of cocommutative
QUE and commutative QFSH algebras. If $U = U(\a)[[\hbar]]$ with
cocommutative coproduct
(where $\a$ is a Lie algebra), then $U'$ is a completion of
$U(\hbar \a[[\hbar]])$; this is a flat deformation of $\wh S(\a)$
equipped with its linear Lie-Poisson structure. If $G$ is a formal group
with function ring $\cO_G$, then $\cO := \cO_G[[\hbar]]$ is a QFSH algebra,
and $\cO^\vee$ is a commutative QUE algebra; it is a quantization of
$(S(\g^*)$, commutative product, cocommutative coproduct, co-Poisson
structure induced by the Lie bracket of $\g)$.

\subsection{Proof that ``twists'' can be made admissible}

\begin{definition}
\label{admissible}
{\it An element $x$ in a QUE algebra
$U$ is admissible if $x\in 1 + \hbar U$, and
if $\hbar \log x$ is in $U' \subset U$.}
\end{definition}

In this subsection, we will prove that for $\gamma, \gamma'$ in
$\Gamma$, the twist $\on{F}_{\gamma,\gamma\gamma'}$ defined in
Proposition \ref{propEH} is twist equivalent to an admissible one.
More precisely, we have
\begin{proposition}
Let $\on{F}_{\gamma,\gamma\gamma'}$ be as Proposition
\ref{propEH}. Then there exists elements
$\on{b}_{\gamma,\gamma\gamma'}$ in $U$ such that
${}^{\on{b}_{\gamma,\gamma\gamma'}}
\on{F}_{\gamma,\gamma\gamma'}:=\on{b}_{\gamma,\gamma\gamma'}^{\otimes
2}
\on{F}_{\gamma,\gamma\gamma'}\Delta_\gamma(\on{b}_{\gamma,\gamma\gamma'}^{-1})$
is admissible.
\end{proposition}
{\em Proof.} Let us denote
$\on{F}_0=\on{F}_{\gamma,\gamma\gamma'}$. We will follow the proof
of Proposition 5.2. in \cite{EH2}: let us construct
$\on{b}=\on{b}_{\gamma,\gamma\gamma'}$ as a product $\cdots
\on{b}_2 \on{b}_1$, where $\on{b}_n\in 1 + \hbar^n U_0$, in such a
way that if $\on{F}_n := {}^{\on{b}_n\cdots \on{b}_1}\on{F}_0$,
then $\hbar\log(\on{F}_n) \in U_0^{\prime\wh\otimes 2} +
\hbar^{n+2} U_0^{\wh\otimes 2}$ (here $U_0$ denotes the
augmentation ideal).

We have already $\hbar\log(\on{F}_0) \in \hbar^2 U_0^{\wh\otimes 2}$.

Expand $\on{F}_0 = 1^{\otimes 2} + \hbar \on{f}_1 + \cdots$, then $\on{Alt}(\on{f}_1) = r$.
Moreover, the coefficient of $\hbar$ in $\on{F}_0^{1,2}\on{F}_0^{12,3} =
\on{F}_0^{2,3}\on{F}_0^{1,23}$ yields $d(\on{f}_1) = 0$, where $d : U(\g)_0^{\otimes 2}
\to U(\g)_0^{\otimes 3}$ is the co-Hochschild differential.
It follows that for some $a_1\in U(\g)_0$, we have $\on{f}_1 = r+d(a_1)$.
Then if we set $\on{b}_1 := \exp(\hbar a_1)$ and $\on{F}_1 = {}^{\on{b}_1}\on{F}_0$,
we get $\on{F}_1 \in 1^{\otimes 2} + \hbar r + \hbar^2 U_0^{\wh\otimes 2}$.
Then $\hbar\log(\on{F}_1) \in \hbar^2 r + \hbar^3 U_0^{\wh\otimes 2}
\subset U_0^{\prime\wh\otimes 2} + \hbar^3 U_0^{\wh\otimes 3}$.

Assume that for $n\geq 2$, we have constructed $\on{b}_1,\ldots,\on{b}_{n-1}$ such that
$\alpha_{n-1} := \hbar\log(\on{F}_{n-1}) \in U_0^{\prime\wh\otimes 2}
+ \hbar^{n+1} U_0^{\wh\otimes 2}$.

Let us recall to technical lemmas from \cite{EH2}:
\begin{lemma} \label{lemma:quot}
The quotient $(U'+\hbar^n U) / (U'+\hbar^{n+1}U)$
identifies with $U(\g)/U(\g)_{\leq n}$. In the same way, the quotient
$(U_0^{\prime\wh\otimes k} + \hbar^n U_0^{\wh\otimes k})
/ (U_0^{\prime\wh\otimes k} + \hbar^{n+1} U_0^{\wh\otimes k})$ identifies with
$U(\g)_0^{\otimes k} / (U(\g)_0^{\otimes k})_{\leq n}$
and the quotient of $\g$-invariant subspaces
$(U_0^{\prime\wh\otimes k} + \hbar^n U_0^{\wh\otimes k})^\g
/ (U_0^{\prime\wh\otimes k} + \hbar^{n+1} U_0^{\wh\otimes k})^\g$
identifies with
$(U(\g)_0^{\otimes k})^\g / (U(\g)_0^{\otimes k})^\g_{\leq n}$.
\end{lemma}
\begin{lemma} \label{lemma:approx}
Assume that $n\geq 2$.
If $f_1,f_2\in (U'_0)^2 + \hbar^{n+1}U_0$ and $g,h\in \hbar^n U_0$,
then $(f_1 + g) \star_\hbar (f_2 + h) = g + h$ modulo
$(U'_0)^2 + \hbar^{n+1}U_0$,
where $\star_\hbar$ is the CBH product for the Lie bracket $[a,b]_\hbar =
[a,b]\slash \hbar$.
\end{lemma}
Let us denote by $\bar\alpha$ the image of the class of
$\alpha_{n-1}$ in $U(\g)_0^{\otimes 2} / (U(\g)_0^{\otimes 2})_{\leq n+1}$
under the isomorphism of this space with $(U_0^{\prime\wh\otimes 2}
+ \hbar^{n+1} U_0^{\wh\otimes 2}) / (U_0^{\prime\wh\otimes 2}
+ \hbar^{n+2} U_0^{\wh\otimes 2})$ (see Lemma \ref{lemma:quot}).
Let $\alpha\in U(\g)_0^{\otimes 2}$ be a representative of $\bar\alpha$, then
$\alpha_{n-1} = \alpha' + \hbar^{n+1}\alpha$, where $\alpha'\in
U_0^{\prime\wh\otimes 2} + \hbar^{n+2} U_0^{\wh\otimes 2}$. Then the twist equation gives
\begin{equation}
\label{equu4}
(-\alpha'-\hbar^{n+1}\alpha)^{1,23} \star_\hbar
(-\alpha'-\hbar^{n+1}\alpha)^{2,3} \star_\hbar
(\alpha'+\hbar^{n+1}\alpha)^{1,2} \star_\hbar
(\alpha'+\hbar^{n+1}\alpha)^{12,3} = 0.
\end{equation}
According to Lemma \ref{lemma:approx}, the image of equality
(\ref{equu4}) in $(U^{\wh\otimes 3} + \hbar^{n+1}
U^{\prime\wh\otimes 3}) / (U^{\wh\otimes 3} + \hbar^{n+2}
U^{\prime\wh\otimes 3}) \simeq U(\g)^{\otimes 3} / (U(\g)^{\otimes
3})_{\leq n+1}$ is $d(\bar\alpha)=0$, where $d$ is the
co-Hochschild differential on the quotient $U(\g)_0^{\otimes
\bullet} / (U(\g)_0^{\otimes\bullet})_{\leq n+1}$. Since $n\geq 2$,
the relevant cohomology group vanishes, so $\bar\alpha =
d(\bar\beta)$, where $\bar\beta\in U(\g)_0 /(U(\g)_0)_{\leq n+1}$.
Let $\beta\in U(\g)_0$ be a representative of $\bar\beta$ and set
$\on{b}_n := \exp(\hbar^n\beta)$, $\on{F}_n :=
{}^{\on{b}_n}\on{F}_{n-1}$, $\alpha_n := \hbar\log(\on{F}_n)$.
Then
$$
\alpha_n = (\hbar^{n+1}\beta)^1 \star_\hbar (\hbar^{n+1}\beta)^2
\star_\hbar \alpha_{n-1} \star_\hbar (-\hbar^{n+1}\beta)^{12}.
$$
According to Lemma \ref{lemma:approx}, the image of $\alpha_n$
in
$$
(U_0^{\wh\otimes 2} + \hbar^{n+1} U_0^{\prime\wh\otimes 2})
/ (U_0^{\wh\otimes 2} + \hbar^{n+2} U_0^{\prime\wh\otimes 2})
\simeq U(\g)_0^{\otimes 2}/(U(\g)^{\otimes 2}_0)_{\leq n+1}
$$
is $\bar\alpha - d(\bar\beta)=0$. So $\alpha_n$ belongs to
$U_0^{\wh\otimes 2} + \hbar^{n+2} U_0^{\prime\wh\otimes 2}$, as
required. This proves the induction step.
\hfill \qed \medskip

\subsection{Proof of Theorem \ref{theoo2}}

Thanks to the previous subsection, we now know that there exists an
element $\on{b}_{\gamma,\gamma\gamma'}$ in $U$ such that ${}^{\on{b}_{\gamma,\gamma\gamma'}}
\on{F}_{\gamma,\gamma\gamma'}:=\on{b}_{\gamma,\gamma\gamma'}^{\otimes 2}
\on{F}_{\gamma,\gamma\gamma'}\Delta_\gamma(\on{b}_{\gamma,\gamma\gamma'}^{-1})$
is admissible. Let us define $$\on{F}_{\gamma,\gamma\gamma'}'={}^{\on{b}_{\gamma,\gamma\gamma'}}
\on{F}_{\gamma,\gamma\gamma'},\hskip1cm \on{i}_{\gamma,\gamma\gamma'}'=\on{i}_{\gamma,\gamma\gamma'} \circ \Ad(\on{b}_{\gamma,\gamma\gamma'}^{-1})$$ and
$$\on{v}_{\gamma,\gamma\gamma',\gamma\gamma'\gamma''}'=
\on{b}_{\gamma,\gamma\gamma'\gamma''}\on{v}_{\gamma,\gamma\gamma',\gamma\gamma'\gamma''}\on{i}^{-1}_{\gamma,\gamma\gamma'}(\on{b}_{\gamma\gamma',\gamma\gamma'\gamma''}^{-1}) \on{b}_{\gamma,\gamma\gamma'}^{-1}.$$
Then it is clear that $\on{F}_{\gamma,\gamma\gamma'}'$,
$\on{i}_{\gamma,\gamma\gamma'}'$ and $\on{v}_{\gamma,\gamma\gamma',\gamma\gamma'\gamma''}'$ still satisfy the conclusion of Theorem \ref{propEH}.

Thanks to the first subsection of this section, applying the
functor {\bf QUE} $\to$ {\bf QFSH} to the algebras
$(U_\gamma,_\gamma,\Delta_\gamma)$ we get algebras
$(U_\gamma',*_\gamma,\Delta_\gamma)$ which are quantizations of
the Poisson algebras $(\cO_{G^*_\gamma},\{-,-\}_\gamma )$. Since
the twists $\on{F}_{\gamma,\gamma\gamma'}'$ are admissible, the
algebra morphisms $\on{i}_{\gamma,\gamma\gamma'}'$ restrict to
the QFSH algebras $U_\gamma'$. Then to end the proof of Theorem
\ref{theoo2}, one has to prove:
\begin{proposition}
\label{prop:v} The elements
$\on{v}_{\gamma,\gamma\gamma',\gamma\gamma'\gamma''}'$ are
admissible.
\end{proposition}
{\em Proof.}
Let us denote $\on{v}=\on{v}_{\gamma,\gamma\gamma',\gamma\gamma'\gamma''}'$.
Suppose $\on{v}$ is not admissible and let $n$ be the bigger
$i$ such that
$\alpha_0 := \hbar\log(\on{v}) \in U_0
+ \hbar^{n+1} U_0$. By the assumption on $v$, we know that $n \geq 2$.
Let us denote by $\bar\alpha$ the image of the class of
$\alpha_0$ in $U(\g)_0/ (U(\g)_0)_{\leq n+1}$
under the isomorphism of this space with $(U_0
+ \hbar^{n+1} U_0) / (U_0
+ \hbar^{n+2} U_0)$ (see Lemma \ref{lemma:quot}).
Let $\alpha\in U(\g)_0$ be a representative of $\bar\alpha$, then
$\alpha_0 = \alpha' + \hbar^{n+1}\alpha$, where $\alpha'\in
U_0 + \hbar^{n+2} U_0$.
Let $f$, $f'$ and $f''$ be respectively the $\hbar \log$ of
$\on{F}_{\gamma,\gamma\gamma'}'$,
$\on{F}_{\gamma\gamma',\gamma\gamma'\gamma''}'$ and
$\on{F}_{\gamma,\gamma\gamma'\gamma''}$. Then the compatibility equation for composition of twists gives
\begin{equation}
\label{equu5}
f''=(\alpha'+\hbar^{n+1}\alpha)^{\otimes 2} \star_\hbar
\on{i}_{\gamma,\gamma\gamma'}^{-1}(f')\star_\hbar f
 \star_\hbar
(-\alpha'-\hbar^{n+1}\alpha)^{12} = 0.
\end{equation}
According to Lemma \ref{lemma:approx}, the image of equality (\ref{equu5})
in $(U^{\wh\otimes 2} + \hbar^{n+1} U^{\prime\wh\otimes 2}) /
(U^{\wh\otimes 2} + \hbar^{n+2} U^{\prime\wh\otimes 2})
\simeq U(\g)^{\otimes 2} / (U(\g)^{\otimes 2})_{\leq n+1}$
is $d(\bar\alpha)=0$. So $\bar\alpha \in \g$ which is a contradiction
with $n\geq 2$.
\hfill \qed \medskip

\section{Example of simple group with action of the Weyl group}

\subsection{Quantization of Majid and Soibelman \cite{MS}}

We start with briefly recalling the Majid and Soibelman's approach
to quantum Weyl group. Let $\g$ be a complex simple Lie algebra,
$U_\hbar(\g)$ be the natural deformation of the universal
enveloping algebra $U(\g)$. Lustig \cite{LU} and Soibelman
\cite{SO} first independently noticed that a simple reflection $w$
in the Weyl group $W$ of $\g$ defines an automorphism $\alpha_w$
on $U_\hbar(\g)$. Then one can extend $U_\hbar(\g)$ by elements
$\bar{w}$ with $\alpha_w(g)=\bar{w}g\bar{w}^{-1}$ for all simple
reflections in $W$. The extended algebra is called by ``quantum
Weyl group" and denoted by $\widetilde{U_\hbar(\g)}$. In \cite{KR}
and \cite{SO}, $\widetilde{U_\hbar(\g)}$ is used to construct
explicit formulas for solutions to the Yang-Baxter equation.

In \cite{MS}, Majid and Soibelman discovered the bicrossed product
structure on $\widetilde{U_\hbar(\g)}$.  Let $w_i$, $1\leq i\leq
\rank(\g)$ be simple reflections in $W$ and $t_j$, $1\leq j\leq
\rank(\g)$ be elements in the maximal torus corresponding to
$\phi_j\left(\begin{array}{cc}-1&0\\ 0&-1\end{array}\right)$ with
$\phi_j:sl_2\hookrightarrow \g$ embedding to the $j$-th vertex of
the Dynkin diagram. Define $\widetilde{W}$ be the group generated
by $w_i$ and $t_j$, which is a covering of the Weyl group $W$ with
the kernel isomorphic to the direct sum of $k$-copies of
$\mathbb{Z}_2$ ($k=\text{rank}(\g)$). The quantum Weyl group
$\widetilde{U_\hbar{(\g)}}$ is proved in \cite{MS}[Corollary 3.4]
to isomorphic to the bicrossed product
\[
k\widetilde{W}^\psi\bowtie_{\alpha, \chi} U_\hbar(\g),
\]
where $\alpha:U_\hbar(\g)\otimes k\widetilde{W}\to U_q(g)$, $\chi:
k\widetilde{W}\otimes k\widetilde{W}\to U_\hbar(\g)$, and
$\psi:k\widetilde{W}\to U_\hbar(\g)\otimes U_\hbar(\g)$ are linear
maps defined by
\[
\alpha(a\otimes wt)=t^{-1}\alpha_w(a)t,\ \ \ \ \ \  \chi(w_1t_1,
w_2t_2)=x^{-1},\ \ \ \ \ \ \psi(wt)=(\bar{w}^{-1}\otimes
\bar{w}^{-1})\Delta\bar{w}.
\]
In the above equation of $\chi$, $x$ is defined to be an element
in $U_\hbar(\g)$ such that
$\alpha_{w_1w_2(\alpha_{w_1}(t_1)t_2)}=\alpha_{w_1t_1}\alpha_{w_2t_2}Ad_{x^{-1}}$
with $x\in U_\hbar(\g)$.

\begin{proposition}\label{prop:ms}
The quantum Weyl group $\widetilde{U_\hbar(\g)}$ is a quantization
of the $\Gamma=\widetilde{W}$ Lie bialgebra $(\g, [\ ,\ ],
\delta)$, where $(\g, [\ ,\ ], \delta)$ is the Lie bialgebra
structure on $\g$ corresponding to the deformation $U_\hbar(\g)$,
and $\widetilde{W}$ acts on $\g$ as the Weyl group ($t$ acts on
$\g$ by adjoint action), and $f_\gamma=\wedge^2(\gamma)\circ
\delta\circ \gamma^{-1}-\delta$ for $\gamma\in \widetilde{W}$.
\end{proposition}
\begin{proof}
Inspired by the above bicrossed product structure on
$\widetilde{U_\hbar(\g)}$, we introduce the following $\Gamma$
quantized universal enveloping algebras for $\Gamma=\widetilde{W}$
generated by the following data.
\begin{itemize}
\item $(U_\hbar(\g)_\gamma, m_\gamma, \Delta_\gamma)=(U_\hbar(\g), m,
\Delta_\gamma)$, where $m$ is the canonical multiplication on
$U_\hbar(\g)$ and $\Delta_\gamma=\alpha(-, \gamma)^{\otimes
2}\circ Ad(\psi(\gamma))\circ \Delta \circ \alpha^{-1}(-, \gamma)$
with $\Delta$ the canonical coproduct on $U_\hbar(\g)$.

\item $i_{\gamma, \gamma\gamma'}:(U_\hbar(\g), m_\gamma)\to
(U_{\hbar}(\g), m_{\gamma\gamma'})$ by $i_{e,
\gamma}=\alpha(-\otimes\gamma):U_\hbar(\g)\to U_\hbar(\g)$ and
$i_{\gamma, \gamma\gamma'}=i_{e, \gamma'}$.

\item $F_{e, \gamma}\in U_\hbar(\g)^{\otimes 2}$ is set to equal
to $\psi(\gamma)$ and $F_{\gamma, \gamma\gamma'}=F_{e,\gamma'}$.
According to \cite{MS}[Lemma 3.3], for any reflection $w_i\in W$,
$F_{e, w_it}=\psi(w_i)=e^{\frac{1}{2}\hbar H_i\otimes
H_i/(\alpha_i, \alpha_i)}(\mathcal{R}_i)^{-1}_{12}=1+\hbar
f_1+O(\hbar^2)$. (Here $(H_i, X_i^+, X_i^-)$ corresponds to the
embedding $\phi_i:sl_2\hookrightarrow\g$ for the $i$-th root
$\alpha_i$ with normal $(\alpha_i, \alpha_i)$.) Because the part
of $e^{\frac{1}{2}\hbar H_i\otimes H_i/(\alpha_i,\alpha_i)}$ is
symmetric, the antisymmetrization of $f_1$ is equal to the
antisymmetrization of the first order term of
$(\mathcal{R}_i)^{-1}_{21}$, which is equal to the definition of
$f_{w_i}$ by the asymptotic expansion of $\mathcal{R}_i$. This
result extends to an arbitrary element $\gamma$ simply because
$w_i$ generates $W$.

\item $v_{e, \gamma, \gamma\gamma'}=\chi(\gamma, \gamma\gamma')\in U_\hbar(\g)^{\otimes
2}$. According to the definition of $\chi(\gamma, \gamma\gamma')$
we see that $v$ can be chosen be an element in
$1+\hbar^2U_\hbar(\g)$ because the $\alpha$ action is associative
up to the $\hbar$-linear terms by \cite{KR}[Formula (13)] and
\cite{KS}[Prop 1.4.10].
\end{itemize}

It is straight forward to check that the cocycle conditions for
$\alpha, \chi, \psi$, and their compatibilities are equivalent to
the conditions for $(U_\hbar, m, \Delta_{\gamma}, i_{\gamma,
\gamma\gamma'}, F_{\gamma, \gamma\gamma}, v_{\gamma,
\gamma\gamma', \gamma\gamma'\gamma''})$ to be a
$\Gamma=\widetilde{W}$ quantized universal enveloping algebra.
Therefore, the corresponding $\Gamma$ quantized universal
enveloping algebra is isomorphic to $\widetilde{U_\hbar(\g)}$.
\end{proof}
\subsection{Admissibility of the twists}
\begin{corollary}
The twists $F_{\gamma, \gamma\gamma'}$ and $v_{\gamma,
\gamma\gamma', \gamma\gamma'\gamma''}$ defined in Proposition
\ref{prop:ms} are admissible. Therefore, the quantum Weyl group
defines a stack of formal series Hopf algebras quantizing the
corresponding stack of Poisson Hopf algebras dual to
$(\widetilde{W}, \g, [\ ,\ ], \delta, f_\gamma)$.
\end{corollary}
\begin{proof}
We look at the formulas for $F_{e, w t}$. According to $\psi$'s
formula, if $w_i$ is a simple reflection, then $F_{e, w_i
t}=e^{\frac{1}{2}\hbar H_i\otimes H_i/(\alpha_i,
\alpha_i)}(\mathcal{R}_i)^{-1}_{12}$. Taking $\hbar\log$ on $F_{e,
w}$, we have
\[
\hbar^2\frac{1}{2}H_i\otimes
H_i/(\alpha_i,\alpha_i)+\hbar\log((\mathcal{R}_i)^{-1}_{12}).
\]

The first term is primitive as $H_i$ is primitive. And the second
term $\hbar\log((\mathcal{R}_i)^{-1}_{12})$ is primitive because
$\hbar\log(\mathcal{R}_i)$ is primitive which was proved in
\cite{EH1}[Theorem 0.1]. Therefore, we conclude that $F_{e,w_it}$
is admissible when $w$ is a simple reflection. And this property
extends to a general element $\gamma$ directly by products.

By Proposition \ref{prop:v}, we also know that $v$ is admissible
because $F$ is admissible.

We conclude the corollary by Theorem \ref{theoo2}.
\end{proof}

\end{document}